\documentclass{conm-p-l}

\begin{document}

\newcommand{\ci}[1]{_{ {}_{\scriptstyle #1}}}

\newcommand{\norm}[1]{\ensuremath{\|#1\|}}
\newcommand{\abs}[1]{\ensuremath{\vert#1\vert}}
\newcommand{\p}{\ensuremath{\partial}}
\newcommand{\pr}{\mathcal{P}}

\newcommand{\pbar}{\ensuremath{\bar{\partial}}}
\newcommand{\db}{\overline\partial}
\newcommand{\D}{\mathbb{D}}
\newcommand{\B}{\mathbb{B}}
\newcommand{\Sp}{\mathbb{S}}
\newcommand{\T}{\mathbb{T}}
\newcommand{\R}{\mathbb{R}}
\newcommand{\C}{\mathbb{C}}
\newcommand{\Z}{\mathbb{Z}}
\newcommand{\N}{\mathbb{N}}
\newcommand{\scrH}{\mathcal{H}}
\newcommand{\scrL}{\mathcal{L}}
\newcommand{\td}{\widetilde\Delta}

\newcommand{\La}{\langle }
\newcommand{\Ra}{\rangle }
\newcommand{\ran}{\operatorname{Ran}}
\newcommand{\Hom}{\operatorname{Hom}}
\newcommand{\tr}{\operatorname{tr}}
\newcommand{\codim}{\operatorname{codim}}
\newcommand{\imag}{\operatorname{Image}}
\newcommand{\vf}{\varphi}
\newcommand{\f}[2]{\ensuremath{\frac{#1}{#2}}}


\newcommand{\entrylabel}[1]{\mbox{#1}\hfill}

\newenvironment{entry}
{\begin{list}{X}%
  {\renewcommand{\makelabel}{\entrylabel}%
      \setlength{\labelwidth}{55pt}%
      \setlength{\leftmargin}{\labelwidth}
      \addtolength{\leftmargin}{\labelsep}%
   }%
}%
{\end{list}}



\numberwithin{equation}{section}

\newtheorem*{thm}{Theorem}
\newtheorem*{prop}{Proposition}
\newtheorem*{cor}{Corollary}

\title{Heat Kernels and Cycles}
\author{Bruno Harris}
\address{Brown University, Department of Mathematics, Box 1917, Providence, RI 02912}
\email{bruno@math.brown.edu}



\maketitle

\section{Introduction}
We use the heat kernel on a compact oriented Riemannian manifold to assign to a sequence of cycles $C_1,\ldots,C_k$ a real number $(C_1,\ldots,C_k)$, provided the $C_i$ satisfy certain conditions.  If $k=2$ then $(C_1,C_2)$ will be the ordinary topological linking number, an integer, while if $k\geq 3$ then $(C_1,\ldots, C_k)$ will involve also iterated integrals of harmonic forms representing the Poincare duals of the cycles $C_i$.  $(C_1,\ldots,C_k)$ will have certain homological invariance properties as functions of the $C_i$.  If the manifold is furthermore K\"ahler then $(C_1,\ldots,C_k)$ has better properties:  it depends only on the complex structure and not on the choice of K\"ahler metric, and the harmonic forms can be replaced by forms which are merely $d$-closed and $d^c$-closed; also, if some of the $C_i$ are complex then $(C_1,\ldots, C_k)$ equals a number $(C_1^*,\ldots,C_i^*), i<k$, in a K\"ahler submanifold.

$(C_1,\ldots,C_k)$ reduced mod $\Z$ can also be defined by using Cheeger--Simons differential characters and their product.  Thus our results using the heat kernel are to obtain an $\R$-valued construction, express it by iterated integrals of harmonic forms, show how passing to homologous cycles introduces specific linking integers, and obtain further results in the K\"ahler case.

In more detail, let $X$ be an $n$-dimensional compact oriented Riemannian manifold with $n$ even, and $K(x,x',t)$ the kernel representing the heat operator 
$$\exp(-t(dd^*+d^*d))$$ 
of the Laplacian $dd^*+d^*d$ on differential forms on $X$.  We regard $K(x,x',t)$ as an $n$-form on $X\times X$ (and function of the parameter $t>0$) such that for each form $\alpha(x)$ on $X$
$$
\int_{x\in X}\alpha(x)\wedge K(x,x',t)=(\exp[-t(dd^*+d^*d)]\alpha)(x').
$$
Since $n=\dim X$ is even, $K(x,x',t)$ is a closed form on $X\times X$ whose cohomology class is Poincare dual to the diagonal $\Delta\subset X\times X$ (see [Harris 1993]).

In our previous work [Harris 1993,2002,2004] we used $K(x,x',t)$ and a related ``linking kernel'' $L(x,x',t)$ do define a real-valued pairing $(A,B)_X$ of cycles $A,\, B$ in $X$ satisfying: $A,\, B$ are disjoint, and
$$
\dim A+\dim B=\dim X-1.
$$
L(x,x',t) satisfies as $n-1$-form on $X\times X$:
\begin{equation}
\label{1.1}
dL(x,x',t)=K(x,x',t)-H(x,x')
\end{equation}
on $X\times X$, where $H(x,x')$ is the kernel for projection to the harmonic forms on $X$.  $X\times X$ is given the product metric and we require $L$ to be orthogonal to harmonic forms on $X\times X$.

The pairing is then defined as the real number
\begin{equation}
\label{1.2}
(A,B)_X = \lim_{t\to 0}\int_{(x,x')\in A\times B} L(x,x',t).
\end{equation}

If both $A,\, B$ bound in $X$ then $(A,B)_X$ is an integer, the ordinary topological linking number.  Thus, if only $A$ bounds then for fixed $A$, $(A,B)_X$ reduced modulo $\Z$ is a function of the homology class of $B$ (with integer coefficients) and defines a function
$$
A\to\Hom(H_*(X;\Z),\R/ \Z)
$$
which is the Abel-Jacobi map.

If $X$ is furthermore complex with K\"ahler metric then we have besides $d$, a ``twisted'' operator
$$
d^c=\f{1}{4\pi}J^{-1}dJ
$$
($J$ is the almost complex structure) and can construct a kernel $l(x,x',t)$ satisfying orthogonality to harmonic forms on $X\times X$ and
\begin{equation}
\label{1.3}
dd^c l(x,x',t)=K(x,x',t)-H(x,x')
\end{equation}
and a pairing
\begin{equation}
\label{1.4}
\La a, b\Ra_X = \lim_{t\to 0}\int_{a\times b}l(x,x',t)
\end{equation}
for complex cycles $a,\,b$ on $X$ with disjoint supports and complex dimensions satisfying
$$
\dim_{\C}a+\dim_{\C}b=\dim_{\C}X-1.
$$
Then $\La a,b\Ra_{X}$ is the Archimedean height pairing.

In [Harris 2002, 2004] we looked at the special case of (\ref{1.2}) where, for a Riemannian manifold $Y$ 
$$
X=Y\times\cdots\times Y=Y^k,\quad k\geq 3,
$$
with product metric on $X$, and the cycles in $X$ are
\begin{eqnarray*}
A & = & \textnormal{the diagonal }\Delta = \{(y,\ldots, y)\}\\
B & = & C_1\times\cdots\times C_k,
\end{eqnarray*}
the $C_i$ being cycles in $Y$ satisfying the following conditions:\\
\newline
\noindent
(1.5)  $C_i=\sum_{r}n_{ir}C_{ir}$ (formal linear combinations with integer coefficients $n_{ir}$) and the $C_{ir}$ are disjoint compact oriented smooth submanifolds of $Y$ of the same codimension $p_i$.  We write: $\abs{C_i}=\textnormal{ support of } C_i=\cup_r C_{ir}=\textnormal{oriented submanifold}$ $\textnormal{of codimension  } p_i$.\\
\newline
\noindent
(1.6) For $i_1,i_2,\ldots,$ distinct, $C_{i_1 r_1},C_{i_2, r_2},\ldots$ intersect transversely.

We distinguish the intersection of cycles, denoted
$$
C_{i_1}\bullet C_{i_2}\bullet\cdots\bullet C_{i_s}
$$
($i_1,i_2,\ldots,i_s$ always distinct), which is a cycle, from the intersection of supports denoted (if $s\geq 2$)
$$
C_{i_1}\cap C_{i_2}\cap\cdots\cap C_{i_s}
$$
meaning $\abs{C_{i_1}}\cap\cdots\cap\abs{C_{i_s}}$.  (The cycle $C_i\bullet C_j$ is defined by the orientation convention that if $C_{ir}, C_{js}$ are oriented by Thom form $v_{ir},v_{js}$ in the normal direction then $v_{ir}\wedge v_{js}$ orients $C_{ir}\cap C_{js}$.)\\
\noindent
The intersection of a cycle $C_{i_1}\bullet\cdots\bullet C_{i_r}$ with a support $C_{j_1}\cap\cdots\cap C_{j_s}$ ($i_1,\ldots,i_r$, $j_1,\ldots,j_s$ distinct indices) is a cycle on this support denoted
$$
C_{i_1}\bullet\cdots\bullet C_{i_r}\bullet (C_{j_1}\cap\cdots\cap C_{j_s}).
$$
>From the disjointness of the components $C_{i1},C_{i2},\ldots,$ of $C_i$ in (1.5) we see that if an intersection, written for simplicity $C_i\bullet\abs{C_j},\ i\neq j$, is homologous to zero (integer coefficients) on $\abs{C_j}$ then $C_i\bullet C_j$ is homologous to zero on $\abs{C_j}$ (the hypothesis means that $C_i=\sum_{u}n_{iu}C_{iu},\ C_j=\sum_{v}n_{jv}C_{jv}$ and $n_{iu}(C_{iu}\bullet C_{jv})$ bounds on $C_{jv}$ for all $v$).\\
\newline
\noindent
(1.7) For $p_i=\textnormal{codimension of } C_i \textnormal{ in } Y$, $\sum_{i=1}^{k}p_i=\dim Y+1=n+1$.

(1.7) together with transversality (1.6) implies that
$$
\abs{C_1}\cap\ldots\cap\abs{C_k}=\textnormal{empty set }\emptyset.
$$
\\
In [Harris 2002, 2004] we further assumed\\
\newline
\noindent
(1.8) For $i=1,\ldots,k,\quad (k\geq 3)$,
$$
\abs{C_1}\cap\ldots\cap\abs{C_{i-1}}\cap\abs{C_{i+1}}\cap\ldots\cap\abs{C_k}=\emptyset
$$
(the intersection of any $k-1$ of the $\abs{C_i}$ is empty).

We now use (\ref{1.2}) to define the pairing $(\Delta,C_1\times\cdots\times C_k)_{Y^k}$ in $Y^{k}$, noting that the diagonal $\Delta$ and $C_1\times\cdots\times C_k$ are disjoint by (1.6) and (1.7) and introduce the notation
\begin{equation}
\tag{1.9}
(C_1,\ldots,C_k)=(\Delta,\,C_1\times\cdots\times C_k)_{Y^k}\quad\textnormal{(assuming (1.6) and (1.7))}.
\end{equation}

We then prove the following theorem in [Harris 2002, 2004]\\
\newline
\noindent 
(1.10) \textbf{ Theorem.} \begin{itshape}Assume (1.5),(1.6),(1.7).  Then
\begin{itemize}
\item[a)] If we interchange $C_r, C_s$ ($r\neq s$) in (1.9) then (1.9) is multiplied by $(-1)^{p_r p_s}$ (where $p_i=\codim C_i$).  This holds even for $k=2$ and does not assume (1.8).
\end{itemize}
\b
For the remaining results we also assume (1.8) and $k\geq 3$.  Then
\begin{itemize}
\item[b)]$(C_1,\ldots,C_k)$ is unchanged if we replace one of the cycles $C_i$ by a homologous cycle $C_i'$ provided $C_1,\ldots,C_{i-1},C_i',C_{i+1},\ldots, C_k$ satisfy the same conditions (1.5)--(1.8) as $C_1,\ldots, C_k$.
\item[c)] Let $\alpha_i$ be the harmonic form Poincare dual (in de Rham cohomology) to the homology class of $C_i$ in $Y$, $i=1,\ldots,k$.  By (1.8) we have 

\item[(1.11)]
$\alpha_1\wedge\ldots\alpha_{i-1}$ restricted to $C_{i+1}\cap\cdots\cap C_k$ is exact:  there exist forms $A_{1\cdots i-1}$ on $C_{i+1}\cap\cdots C_k$ such that $\alpha_1\wedge\ldots\wedge\alpha_{i-1}=d A_{1\cdots i-1}$, for $i=2,\ldots, k$.  For $i=k$ we interpret this to mean:  $\alpha_1\wedge\ldots\wedge\alpha_{k-1} =d A_{1\cdots k-1}$ on $Y$.

We then have:

\item[(1.12)] For any choice of forms $A_{1\cdots i-1}$ on $C_{i+1}\cap\cdots\cap C_k$ as above,
\begin{equation}
\tag{1.13}
(C_1,\ldots, C_k)  =  (-1)^{q+1}[\int_{C_3\bullet\cdots\bullet C_k}(A_1\alpha_2-A_{12})
\end{equation}
\begin{equation*}
 +\int_{C_4\bullet\cdots\bullet C_k}(A_{12}\alpha_3-A_{123})+\cdots+\int_{Y}A_{1\cdots k-1}\alpha_k]
\end{equation*}
where $q=\sum_{r<s}p_rp_s$.  (If we choose any $n$-form $A_{1\cdots k}$ on $Y$ such that $\int_{Y}A_{1\cdots k}=0$, so $d A_{1\cdots k}=\alpha_1\wedge\cdots\alpha_k=0$ for degree reasons, and let $C_{k+1}=Y$, we can write
$$
(C_1,\ldots C_k)=(-1)^{q+1}\sum_{i=2}^{k}\int_{C_{i+1}\bullet\cdots\bullet C_{k+1}}(A_{1\cdots i-1}\alpha_i-A_{1\cdots i})).
$$
As example, for $k=3$, 
$$(C_1, C_2, C_3)=(-1)^{q+1}\left[\int_{C_3}(A_1\alpha_2-A_{12})+\int_Y A_{12}\alpha_3\right]$$.
\item[d)] If $Y$ is complex and the metric is K\"ahler then $(C_1,\ldots, C_k)$ depends only on the complex structure and not on the choice of K\"ahler metric.
\end{itemize}
\end{itshape}
Remarks:  b) follows from formula (1.13) in c) and from a) since (1.13) shows that $(C_1,\ldots, C_k)$ does not depend on $C_1$ or $C_2$ but depends only on $\alpha_1,\ldots, \alpha_k$ and $C_3,\ldots, C_k$.  The homology statment b) was stated in a possibly misleading way in [Harris 2004, line 11b of Th. 3.2]

In the next section we will obtain an improvement of the homology invariance statement b) for $k=3$ or $4$ by generalizing the disjointness condition of (1.8) to allow intersections homologous to zero, and bringing in certain integer linking numbers.  We also obtain improvement of the K\"ahler case d), where we no longer need harmonic forms -- a large practical advantage.

\section{}
We now state our more general assumptions and the corresponding results.  The notation $Y, C_1,\ldots, C_k$ is as in section 1, and the cycles $C_1,\ldots C_k$ will always be assumed to satisfy (1.5), (1.6) (transversality), and (1.7) (sum of codimensions is $\dim Y+1$).  Then the intersection of supports $C_1\cap\cdots\cap C_k$ is empty and so $(C_1,\ldots, C_k)$ is defined (even for $k=2$; usually we assume $k\geq 3$) and we claimed that it satisfies skew-symmetry; however we need the proof of this for our next results so we will sketch it now:\\

\noindent
(2.1) \textbf{ Proposition.}\begin{itshape} Assuming (1.5), (1.6), (1.7) and $k\geq 2$, $(C_1,\ldots, C_k)$ is multiplied by $(-1)^{p_r p_s}$ if $C_r,C_s$ ($r\neq s$), of codimensions $p_r, p_s$, are interchanged.
\end{itshape}
\begin{proof}  (see [Harris, 2004], p.89).  In the definition of (\ref{1.2}) we may use any kernel $L=L(x,x',t)$ satisfying (\ref{1.1}) and orthogonality to harmonic forms on $X\times X$, since any two such $L$ differ by an exact form.  Since $K(x,x',t)$ on $X\times X=Y^k\times Y^k$ is the product of the corresponding $K_i(y,y',t)$ on $Y\times Y$ where $x=(y_1,\ldots,y_k)$, $x'=(y_1',\ldots,y_k')$ and $K_i=K(y_i,y_i',t)$, so that $K=K_1\cdot\ldots\cdot K_k$ and similarly $H=H_1\cdot\ldots\cdot H_k$, ($H_i$ on $Y\times Y$), we may choose
\begin{equation}
\tag{2.2}
L(x,x',t)=L_1K_2\cdot\ldots\cdot K_k+H_1L_2K_3\cdot\ldots\cdot K_k+\cdots+H_1\cdot\ldots\cdot H_{k-1}L_k,
\end{equation}
where $L_i=L_Y(y_i,y_i',t)$ and $d L_{Y}=K_Y-H_Y$ on $Y$.  However we may choose a different $L$ on $Y^{k}\times Y^{k}$, for instance
$$
L'(x,x',t)=K_1L_2K_3\cdot\ldots\cdot K_k+L_1H_2K_3\cdot\ldots\cdot K_k+\sum_{i=3}^{k}H_1\cdot\ldots\cdot H_{i-1}L_i K_{i+1}\cdot\ldots\cdot K_{k}.
$$
Next we have to integrate $L$ or $L'$ over $\Delta\times (C_1\times\cdots\times C_k)$ and take limit as $t\to 0$.  Comparison of these two integrals then gives the factor $(-1)^{p_rp_s}$ for interchanging $C_1, C_2$, ending the proof of (2.1).
\end{proof}
However we will need more of the proof (loc. cit. p. 86-87).  We first integrate $L$ in the second variable over $C_1\times\ldots\times C_k$.  The integration of $L_Y(y,y',t)$ over $y'\in C_i$ gives a form $\Gamma_{C_i,t}(y)$ satisfying , on $Y$, 
$$
d\Gamma_{C_i,t}=K_{C_i,t}-\alpha_i
$$
where $K_{C_i,t}$ is the integral of $K_Y(y,y',t)$ over $y'\in C_i$ and $\alpha_i$ is the integral of $H_Y(y,y')$ over $y'\in C_i$.  Thus the integral of $L(x,x',t)$ over $x'\in C_1\times\cdots\times\ C_k$ and then restriction to $x=(y_1,\ldots,y_k)\in\Delta$ give the following form on $\Delta=Y$:
\begin{equation}
\tag{2.3}
(-1)^q\sum_{i=1}^{k}(-1)^{p_1+\cdots+p_{i-1}}\alpha_1\cdots\alpha_{i-1}\Gamma_{C_i,t}K_{C_{i+1},t}\cdot\ldots\cdot K_{C_k,t}\quad (q=\sum_{r<s}p_rp_s).
\end{equation}
$\alpha_i$ is again the harmonic form Poincare dual to $C_i$.  The forms $K_{C_i,t}$ and $\Gamma_{C_i,t}$ have the following behavior for $t\to 0$ (obtained from the heat kernel asymptotics of $K_Y(y,y',t)$):  as $t\to 0$ and for $y$ outside $C_i$, $K_{C_i,t}(y)$ approaches $0$ ``exponentially'', namely like $e^{-\f{r^{2}}{4t}}$, $r=\textnormal{distance to } C_i$, whereas on $C_i$ both $K_{C_i}$ and $\Gamma_{C_i}$ approach $\infty$ like a negative power of $t$.  $K_{C_i,t}$ regarded as a current approaches the Dirac current of integration over $C_i$, namely $\delta_{C_i}$, as $t\to 0$.  $\Gamma_{C_i,t}$ approaches a current which is also singular on $C_i$ and near $C_i$ is the normalized angular form measuring  ``angular measure'' on small spheres normal to $C_i$ (like $\f{1}{2\pi}d\theta$ around the origin in $\R^2$).  Outside $C_i$, $\Gamma_{C_i}$ is smooth.

One choice for $L(x,x',t)$ on $X\times X$ is 
\begin{equation}
\tag{2.4}
L(x,x',t)=\f{1}{2}\int_{\tau=t}^{\infty}d^{*}K(x,x',\tau)d\tau
\end{equation}
which is in the image of $d^*$ on $X\times X$, and similarly for $L_Y(y,y',t)$ on $Y\times Y$.  We can also show that $\Gamma_{C_i,t}(y)$, obtained by integrating $L_Y$ over $y'\in C_i$, is in the image of $d^*$ on $Y$:  we start with
\begin{equation}
\tag{2.5}
K(y,y',t)=\sum_{\lambda}e^{-\lambda t}(\ast\varphi_{\lambda}(y))\wedge\varphi_{\lambda}(y')
\end{equation}
$\lambda =$ eigenvalues of the Laplacian and $\varphi_\lambda$ the corresponding eigenforms, $\ast=$ Hodge star operator.  Next we apply $d^*$ on $Y\times Y$ to (2.5) and note that for the product metric on $Y\times Y$, $d^*$ obeys the Leibniz product rule for products of the form $\varphi(y)\wedge\psi(y')$.  Thus $d^* K(y,y',t)$ is an infinite series of terms $e^{-\lambda t}d^* \psi(y)\wedge\mu(y')$, and if we integrate over $y'\in C_i$ the result is in image of $d^*$ on $Y$.  Thus we have 
\begin{equation}
\tag{2.6} 
\Gamma_{C_i,t}(y) \textnormal{ obtained from (2.4) is in image of } d^* \textnormal{ on } Y.
\end{equation}

We use this discussion, involving (2.2) and (2.3), for the case $k=2$ where $C_1, C_2$ are disjoint cycles on $Y$ with codimensions $p_1,p_2$ satisfying $p_1+p_2=n+1=\textnormal{odd integer}$ so that $(-1)^{p_1p_2}=1$.  Then by (2.3) we obtain for $(C_1,C_2)=(\Delta, C_1\times C_2)_{Y\times Y}$
\begin{equation}
\tag{2.7} 
(C_1, C_2)=\lim_{t\to 0}\int_{Y}[\Gamma_{C_1,t}(t)K_{C_2,t}(y)+(-1)^{p_1}\alpha_1\wedge\Gamma_{C_2,t}].
\end{equation}
The integral of the first term over $Y$ approaches $\int_{C_2}\Gamma_{C_1}$, since $K_{C_2,t}$ approaches the $\delta$-function $\delta_{C_2}$ as $t\to 0$, and $\Gamma_{C_1}$ is continuous on a neighborhood of $C_2$ as $C_1\cap C_2=\emptyset$.  The integral of the second term over $Y$ is zero since $\Gamma_{C_2,t}$ is orthogonal to harmonic forms.  Thus we have (using also skew-symmetry)
\begin{equation}
\tag{2.8}
(C_1,C_2)=\lim_{t\to 0}\int_{C_2}\Gamma_{C_1,t}=\int_{C_2}\Gamma_{C_1}=\int_{C_1}\Gamma_{C_2}=(C_2,C_1).
\end{equation}

Suppose further that $C_1, C_2$ both bound in $Y$, say $C_i=\p D_i$, with $D_1$ transverse to $C_2$ and $D_2$ transverse to $C_1$.  Then $\alpha_1, \alpha_2$ are both $0$ and $d\Gamma_{C_i,t}=K_{C_i,t}$ for $i=1,2$ (and $t>0$).  Now (2.8) says $(C_1, C_2)=\lim_{t\to 0}\int_{D_2}K_{C_1,t}=$`` intersection number'' of $C_1$ with interior of $D_2$, so we may use the terminology
\begin{equation}
\tag{2.9}
(C_1,C_2)= \textnormal{Linking number } Lk(C_1,C_2)_Y=Lk(C_2,C_1)\in\Z.
\end{equation}

We can also see this linking from the angular behavior of $\Gamma_{C_i}$ around $C_i$.  We can now return to the calculation of $(C_1,\ldots, C_k)$ for $k\geq 3$ under the following condition (2.10) which generalizes (1.8).  (1.8) was that any $k-1$ of the $C_i$ have empty intersection.

For cycles $Z,W$ with transverse supports $\abs{Z},\abs{W}$, denote $Z\bullet\abs{W}$ the intersection of the cycle $Z$ with the support of $W$, which is a cycle on $\abs{W}$.  Thus, as discussed after (1.6) in Section 1
$$
C_1\bullet\cdots\bullet C_{i-1}\bullet\abs{C_{i+1}\cap\ldots\cap C_k}
$$
denotes the intersection of the cycle $C_1\bullet\cdots\bullet C_{i-1}$ with the support $\abs{C_{i+1}}\cap\ldots\cap\abs{C_k}$ for $i=2,\ldots, k-1$.  For $i=k$ this notation means just the cycle $C_1\bullet\cdots\bullet C_{k-1}$ on $Y$.

\begin{itemize}
\item[(2.10)] Assumption.  The $C_i$ satisfy (1.5)-(1.7) and:
$$
C_2\bullet (C_3\cap\cdots\cap C_k)
$$
bounds on $C_3\cap\cdots\cap C_k$.  For $i=2,\ldots, k$
$$
C_1\bullet\cdots\bullet C_{i-1}\bullet(C_{i+1}\cap\cdots\cap C_k)
$$
bounds on $C_{i+1}\cap\cdots\cap C_k$ (for $i=k$:  $C_1\bullet\cdots\bullet C_{k-1}$ bounds on $Y$).
\end{itemize}
In general, for $Z$, $W$ as above, the condition
\begin{itemize}
\item[(2.11)] $Z,W$ meet transversely and $Z\bullet\abs{W}$ bounds in $\abs{W}$
\end{itemize}
is equivalent to
\begin{itemize}
\item[(2.12)]Any integer cocycle $\xi$ representing the Poincare dual to the homology class of $Z$ (integer coefficients), cobounds when restricted to the support of $\abs{W}$: $\xi= \delta A$, $A=$ cochain on $\abs{W}$.  (Note that $\xi$ restricted to $\abs{W}$ is Poincare dual to the homology class of $Z\bullet\abs{W}$).
\end{itemize}

Another equivalent condition is that $Z$ is homologous on $X$ to a cycle $Z'$ on $X\setminus\abs{W}$ (using Poincare-Lefschetz duality).

Thus (2.10) implies that 
\begin{itemize}
\item[(2.13)] $\alpha_1\wedge\cdots\wedge\alpha_{i-1}$ restricted to $C_{i+1}\cap\cdots\cap C_k$ is exact:
$$
\alpha_1\wedge\cdots\alpha_{i-1}=d A_{1\ldots i-1}\ \textnormal{on this set.}
$$
For $i=k$, $\alpha_1\wedge\cdots\wedge\alpha_{k-1}$ is exact on $Y$.
\end{itemize}
In (2.10) we assumed $C_2\bullet(C_3\cap\cdots\cap C_k)$ bounds on $C_3\cap\cdots\cap C_k$.  By the discussion following (1.6) we see that this implies the cycle $C_2\bullet C_3\bullet\cdots\bullet C_k$ bounds on $C_3\cap\cdots\cap C_k$.  Also (2.10) assumes that $C_1\bullet (C_3\cap\cdots\cap C_k)$ bounds on $C_3\cap\cdots\cap C_k$.  Further, these two cycles on $C_3\cap\cdots\cap C_k$ are disjoint as $C_1\cap C_2\cap\cdots\cap C_k=\emptyset$.  Thus we have a well defined linking number (integer) (as in (2.9))
$$
Lk(C_1\bullet(C_3\cap\cdots\cap C_k),C_2\bullet C_3\bullet\cdots\bullet C_k)_{C_3\cap\cdots\cap C_k}
$$
which we will abbreviate as 
$$
Lk(C_1, C_2\bullet C_3\bullet\cdots\bullet C_k)_{C_3\cap\cdots\cap C_k}.
$$

We can now state\\
\noindent
(2.14) \textbf{ Proposition.}\begin{itshape} Suppose $C_1,\ldots,C_k$ on $Y$ satisfy (2.10) (including (1.5)-(1.7) and $k\geq 3$).  Then
\begin{equation}
\tag{2.15}
(C_1,\ldots,C_k)=(-1)^q[Lk(C_1,C_2\bullet C_3\bullet\cdots\bullet C_k)_{C_3\cap\cdots\cap C_k}-
\end{equation}
\begin{equation*}
\sum_{i=2}^{k}\int_{C_{i+1}\bullet\cdots\bullet C_{k+1}}(A_{1\ldots i-1}\alpha_i-A_{1\ldots i})]
\end{equation*}
where for $i\leq k$, $A_{1\ldots i-1}$ are any forms satisfying (2.13), $A_{1\ldots k}$ is any form on $Y$ satisfying $\int_{Y}A_{1\ldots k}=0$, and $C_{k+1}=Y$, for $i=k$, $C_{i+1}\bullet\cdots\bullet C_{k+1}=C_{k+1}=Y$.
\end{itshape}

\begin{proof}
We start with the sum in (2.3),
$$
\sum_{i=1}^{k}(-1)^{p_1+\cdots+p_{i-1}}\alpha_1\cdots\alpha_{i-1}\Gamma_{C_i,t}K_{C_{i+1},t}\cdots K_{C_k,t}
$$
which is integrated over $Y$, and the limit $t\to 0$ is taken.  The first term gives
$$
\int_{Y}\Gamma_{C_1,t}K_{C_2,t}\cdots K_{C_k,t}.
$$
In view of the discussion following (2.12) with $Z=C_1$ and $\abs{W}=C_3\cap\cdots\cap C_k$, we know that $Z$ is homologous to a cycle $Z'$ not meeting $C_3\cap\cdots\cap C_k$ and so the Poincare dual form $\alpha_1$ to $Z'$ is exact outside $Z'$, thus there is a neighborhood $U$ of $C_3\cap\cdots\cap C_k$ and a form $A_1$ on $U$ such that
$$
\alpha_1=d A_1\ \textnormal{on } U.
$$
Since $K_{C_2,t}\cdots K_{C_k,t}\to 0$ rapidly outside $C_2\cap\cdots\cap C_k$ as $t\to 0$, and $C_2\cap\cdots\cap C_k\subset U$, we may neglect $\int_{Y\setminus U}\Gamma_{C_1,t}K_{C_2,t}\cdots K_{C_k,t}$.  For $t$ near $0$, $\int_{U}\Gamma_{C_1,t}K_{C_2,t}\cdots K_{C_k,t}$ is close to $\int_{C_2\bullet\cdots\bullet C_k}\Gamma_{C_1,t}$ (noting $C_1$ is disjoint from $C_2\cap\cdots\cap C_k$).  Since by (2.10) $C_2\bullet\cdots\bullet C_k=\p D$, $D$ a chain on $C_3\cap\cdots\cap C_k$, $\int_{C_2\bullet\cdots\bullet C_k}\Gamma_{C_1,t}=\int_{D}(K_{C_1,t}-\alpha_1)$ and $\int_{D}K_{C_1,t}$ approaches the linking number
$$
Lk(C_1, C_2\bullet C_3\bullet\cdots\bullet C_k)_{C_3\cap\cdots\cap C_k}
$$
as $t\to 0$.  Thus the first integral approaches (this linking number ) $-\int_{C_2\bullet\cdots\bullet C_k}A_1$.

We proceed similarly for the next term
$$
(-1)^{p_1}\int_{Y}\alpha_1\Gamma_{C_2,t}K_{C_3,t}\cdots K_{C_k,t}.
$$
With the same neighborhood $U$ of $C_3\cap\cdots\cap C_k$ and form $A_1$ as above, we consider $U$ as manifold with boundary $\p U$ and see that the integral over $Y\setminus U$ approaches $0$ with $t$.  On $U$ the integral is 
\begin{eqnarray*}
(-1)^{p_1}\int_U dA_1\Gamma_{C_2,t}K_{C_3,t}\cdots K_{C_k,t} & = & (-1)^{p_1}\int_{\p U}A_1\Gamma_{C_2,t} K_{C_3,t}\cdots K_{C_k,t}\\
 &  & +\int_{U} A_1d\Gamma_{C_2,t}K_{C_3,t}\cdots K_{C_k,t}.
\end{eqnarray*}
The first integral $\to 0$ with $t$, and the second, in which we have $d\Gamma_{C_2,t}=K_{C_2,t}-\alpha_{C_2}$ approaches
$$
\int_{C_2\bullet\cdots\bullet C_k} A_1-\int_{C_3\bullet\cdots\bullet C_k} A_1\alpha_2.
$$
The sum of the first two integrals is thus
$$
Lk(C_1,C_2\bullet\cdots\bullet C_k)-\int_{C_3\bullet\cdots\bullet C_k}A_1\alpha_2.
$$

In the same way we evaluate the remaining integrals and obtain the right hand side of (2.15), with a particular choice of $A_1, A_{12},\ldots$ each defined on a neighborhood of $C_3\cap\cdots\cap C_k, C_4\cap\cdots\cap C_k,\ldots,C_k$.  To see that (2.15) is unchanged if we replace each $A_{1\ldots i-1}$ by an $A'_{1\ldots i-1}$ defined on $C_{i+1}\cap\cdots\cap C_k$ (rather than on a neighborhood), write
$$
A'_{1\ldots i-1}=A_{1\ldots i-1}+B_{i-1}
$$
$B_{i-1}$ a closed form on $C_{i+1}\cap\cdots\cap C_k$, and note that $\int_{C_3\bullet\cdots\bullet C_k}B_1\alpha_2=0$ because $B_1$ is closed and $\alpha_2$ is exact on $C_3\cap\cdots\cap C_k$, also $\int_{C_3\bullet\cdots\bullet C_k}B_2-\int_{C_4\bullet\cdots\bullet C_k}B_2\alpha_3=0$ because $B_2$ is closed and $\alpha_3$ is the Poincare dual form to $C_3$ in $X$ and hence to
$$
C_3\bullet(C_4\cap\cdots\cap C_k)\ \textnormal{in } C_4\cap\cdots\cap C_k.
$$
Finally all the integrals with $B_i$'s add up to zero, proving Proposition (2.14)
\end{proof}
\noindent
(2.16)\textbf{ Corollary.}\begin{itshape} Suppose $C_1,\ldots, C_k$ satisfy the conditions of Proposition (2.14).  Suppose $C_1', C_2'$ are cycles homologous to $C_1, C_2$ respectively and $C_1', C_2', C_3, \ldots, C_k$ satisfy the general conditions (1.5)-(1.7) (basically the transversality).  Then $C_1', C_2',$ $C_3,\ldots, C_k$ also satisfy (2.10) and therefore satisfy (2.14).  Further,
\begin{eqnarray*}
\lefteqn{(C_1',C_2',C_3,\ldots, C_k) - (C_1,C_2,C_3,\ldots, C_k)=} & &\\
 & &
  -(-1)^q[Lk(C_1,C_2\bullet C_3\bullet\cdots\bullet C_k)_{C_3\cap\cdots\cap C_k}-Lk(C_1',C_2'\bullet C_3\bullet\cdots\bullet C_k)_{C_3\cap\cdots\cap C_k}]. 
\end{eqnarray*}
\end{itshape}
\begin{proof}
The terms on the right hand side of (2.15) involve in the integrals only the $\alpha_i$ and $A_i$ forms on $C_j\cap\cdots\cap C_{k+1}$ for $j\geq 3$ and so are unchanged if $C_1, C_2$ are replaced by $C_1',C_2'$.  Condition (2.10) for $C_1,\ldots, C_k$ implies (2.10) for $C_1',C_2',C_3,\ldots,C_k$ because $C_1',C_2'$ are homologous to $C_1, C_2$, concluding the proof.
\end{proof}

We will introduce now a condition like (2.10) but independent of the order of the indices $1,\ldots, k$.
\begin{itemize}
\item[(2.17)] Assumption.  $C_1,\ldots, C_k$ satisfy (1.5)-(1.7) and:  if $i_1,\ldots, i_r,j_1,\ldots, j_s$ are $k-1$ distinct indices between $1$ and $k$, with $r\geq 1$, and $s\geq 0$, then each intersection cycle $C_{i_1}\bullet\cdots\bullet C_{i_r}\bullet (C_{j_1}\cap\cdots\cap C_{j_s})$ is homologous to zero on $C_{j_1}\cap\cdots\cap C_{j_s}$.  If $r=k-1,s=0$ this means $C_{i_1}\bullet\cdots\bullet C_{i_{k-1}}$ is homologous to zero on $Y$.
\end{itemize}

Note that (2.17) implies (2.10).  As example, $C_2\bullet (C_3\cap\cdots\cap C_k)$ bounds on $C_3\cap\cdots\cap C_k$ and so $C_2\bullet C_3\bullet\cdots\bullet C_k$ also bounds on $C_3\cap\cdots\cap C_k$.

Further, if $C_1',\ldots, C_{i-1}'$ are homologous to $C_1,\ldots, C_{i-1}$ respectively and if the cycles $C_1',\ldots, C_{i-1}',C_i,\ldots, C_k$ satisfy (1.5)-(1.7) then if $C_1,\ldots, C_k$ satisfy (2.17), $C_1'\bullet\cdots\bullet C_{i-1}'\bullet(C_{i+1}\cap\cdots\cap C_k)$ is homologous to zero on $C_{i+1}\cap\cdots\cap C_k$.

Our result is now\\
\newline
\noindent
(2.18)\textbf{ Theorem.}\begin{itshape}
Let $C_1,\ldots, C_k$ ($k\geq 3$), be cycles in $Y$, a compact oriented Riemannian manifold of even dimension $n$, satisfying (1.5)-(1.7) and (2.17).  Then the expression $(C_1,\ldots, C_k)$ satisfies
\begin{itemize}
\item[a)] ``Skew-symmetry'':  if $C_i$ and $C_j$ are interchanged for $i\neq j$ then $(C_1,\ldots, C_k)$ is multiplied by $(-1)^{p_i p_j}$, $p_i=\codim C_i$.
\item[b)] Formula (2.15), a particular integer linking number added to a sum of iterated integrals involving the harmonic forms $\alpha_i$ Poincare dual to the $C_i$.
\item[c)] Suppose $C_1',\ldots,C_k'$ also satisfy (1.5)-(1.7) and (2.17) and $C_i'$ is homologous to $C_i$ for $i=1,\ldots,k$.  Suppose also that for every set of distinct indices $i_1,\ldots, i_r$ between $1$ and $k$, each intersection of supports
$$
\abs{C_{i_1}}\cap\cdots\cap\abs{C_{i_m}}\cap\abs{C_{i_{m+1}}'}\cap\cdots\cap\abs{C_{i_r}'}
$$
is transverse (``mutual transversality'' of the $C_i, C_j'$).  If $k=3$ or $4$, then $(C_1,\ldots, C_k)=(C_1',\ldots,C_k')+$ a specific sum of integer linking numbers involving the $C_i$ and $C_j
'$.
\item[d)]  Suppose $Y$ is furthermore complex and the metric is K\"ahler.  Then in (2.15) in the integrals on the right hand side we may replace each real harmonic form $\alpha_i$ by any de Rham cohomologous form $\alpha_i'$ satisfying $d\alpha_i'=0=d^c\alpha_i'$ and replace the $A_{1\ldots i}$ by any forms $A'_{1\ldots i}$ satisfying (2.13) with respect to the $\alpha_i'$.  Thus $(C_1,\ldots, C_k)$ is independent of the choice of K\"ahler metric and depends only on the complex structure of $\ Y$.  Furthermore if for some $i\leq k$, the supports $C_i\cap\cdots\cap C_k, C_{i-1}\cap\cdots\cap C_k,\ldots, C_k$ are all complex submanifolds of $Y$ then the integrals over them (and over $Y$) are zero and the whole construction of $(C_1,\ldots, C_k)$ in $Y$ reduces to the lower dimensional expression $(C_1^*,C_2^*,\ldots, C_{i-1}^*)$ in $Y^*$ where $C_j^*=C_j\bullet C_i\bullet\cdots\bullet C_k$ and $Y^*=C_i\cap\cdots\cap C_k$.
\end{itemize} 
\end{itshape}
\begin{proof}
Part a) has already been proved.  Further, assumption (2.17) implies (2.10) and so also implies (2.13)-(2.16).  Thus we only have to prove c) and d).

In c), suppose $k=3$ and $C_1,C_2,C_3$ and $C_1',C_2',C_3'$ satisfy the conditions of c).  Then by (2.16),
\begin{eqnarray*}
(C_1,C_2,C_3) & = & (C_1',C_2',C_3)+(-1)^q[Lk(C_1,C_2\bullet C_3)_{C_3}-Lk(C_1',C_2'\bullet C_3)_{C_3}]
\end{eqnarray*}
\begin{eqnarray*}
(C_1',C_2',C_3) & = & (-1)^{p_1 p_3}(C_3,C_2',C_1')\ \textnormal{by a)}\\
 & = & (-1)^{p_1 p_3}[(C_3',C_2',C_1')+(-1)^q[Lk(C_3,C_2'\bullet C_1')_{C_1'}-(-1)^qLk(C_3',C_2'\bullet C_1')_{C_1'}]\\
  & = & (C_1',C_2',C_3')+(-1)^{q+p_1p_3}[Lk(C_3,C_2'\bullet C_1')_{C_1'}-Lk(C_3',C_2'\bullet C_1')_{C_1'}]
\end{eqnarray*}
Finally,
\begin{eqnarray*}
(C_1,C_2,C_3)-(C_1',C_2',C_3')& = & (-1)^q[Lk(C_1,C_2\bullet C_3)_{C_3}-Lk(C_1',C_2'\bullet C_3)_{C_3}\\
 &  & +(-1)^{p_1p_3}(Lk(C_3,C_2'\bullet C_1')_{C_1'}-Lk(C_3',C_2'\bullet C_1')_{C_1'})].
\end{eqnarray*}

For $k=4$, we proceed similarly:  $(C_1,C_2, C_3, C_4)=(C_1',C_2',C_3,C_4)$ modulo certain linking numbers given in (2.16).  Similarly $(C_1', C_2',C_3',C_4')-(C_1',C_2',C_3,C_4)$ is a sum of linking numbers, by (2.16) and skew-symmetry.

This prove c) for $k=3,4$ but the proof does not extend to $k>4$.

Proof of d):  Let $Y$ be K\"ahler and choose any real forms $\alpha_i'$ which are both $d$ and $d^c$ closed and represent Poincare dual cohomology classes to the homology classes of the $C_i$, recalling that
$$
(\ker d)\cap(\ker d^c)/\imag dd^c\to\ker d/\imag d
$$
is an isomorphism on a K\"ahler manifold.

Let $A'_{1\ldots i-1}$ be forms satisfying (2.13) for the $\alpha_i'$:
$$
\alpha_1'\wedge\cdots\wedge\alpha_{i-1}'= dA'_{1\ldots i-1}\ \textnormal{on } C_{i+1}\cap\cdots\cap C_k
$$
and $\alpha_1'\wedge\cdots\wedge\alpha_{k-1}'= dA'_{1\ldots k-1}$ on $X$.  The $A'$ forms exist by (2.10), the same reason that the $A$'s exist.  Let $A'_{1\ldots k}$ satisfy $\int_{Y}A'_{1\ldots k}=0$.

The previously defined harmonic forms $\alpha_i$ also satisfy $d\alpha_i=0, d^c\alpha_i=0$, and $\alpha_i-\alpha_i'$ is assumed to be $d$-exact.  Thus the $d,d^c$ lemma applies to the $d^c$ closed and $d$-exact forms $\alpha_i-\alpha_i'$ and says that there exist forms $\beta_i$ on $Y$ such that
$$
\alpha_i=\alpha_i'+dd^c\beta_i,\quad i=1,\ldots,k.
$$

Having chosen the forms $A'$ as above, let 
\begin{equation}
\tag{2.19}
A_{1\ldots i-1}=A'_{1\ldots i-1}+d^c\sum_{j=1}^{i-1}\alpha_1\cdots\alpha_{j-1}\beta_j\alpha_{j+1}'\cdots\alpha_{i-1}'
\end{equation}
(on the same submanifold).  Then the $A$ forms satisfy (2.13).  In particular, 
$$\int_Y A_{1\ldots k}=\int_{Y}A'_{1\ldots k}=0,$$ 
since the integral over $Y$ of a form $d^c\gamma$ is zero.

>From (2.19) for $i, i+1$, we get
\begin{equation}
\tag{2.20}
(A_{1\ldots i}-A_{1\ldots i-1}\alpha_i)-(A'_{1\ldots i}-A'_{1\ldots i-1}\alpha_i')=d[(-1)^{p_1+\cdots p_{i-1}}A'_{1\ldots i-1}\wedge d^c\beta_i]
\end{equation}
and so both sides of (2.20) have integral equal to zero over the cycle $C_{i+1}\bullet\cdots\bullet C_{k+1}$ (where $C_{k+1}=Y$).  This shows that (2.15) is satisfied by the $\alpha_i'$ and $A_i'$ as well as by the harmonic $\alpha_i$ and the $A_i$.  Henceforth for simplicity we write $\alpha_i, A_{1\ldots i}$ instead of $\alpha_i', A'_{1\ldots i}$.

To prove the last part of d), assume that $C_{i+1}\cap\cdots\cap C_{k+1}$ (i.e. $Y$ for $i=k$), is a complex submanifold.  Then $\alpha_1\wedge\cdots\wedge\alpha_{i-1}$ restricted to this submanifold is $d$-exact and $d^c$-closed and so is $dd^c$ exact:  we may choose $A_{1\ldots i-1}=d^c B_{1\ldots i-1}$ on $C_{i+1}\cap\cdots\cap C_{k+1}$.  Thus $A_{1\ldots i-1}\wedge\alpha_i= d^c(B_{1\ldots i-1}\wedge\alpha_i)$ and $\int_{C_{i+1}\bullet\cdots\bullet C_{k+1}}A_{1\ldots i-1}\alpha_i=0$.  If furthermore $C_{i+2}\cap\cdots\cap C_{k+1}$ is also complex, we may choose $A_{1\ldots i}=d^c B_{1\ldots i}$ and so
$$
\int_{C_{i+1}\bullet\cdots\bullet C_{k+1}}(A_{1\ldots i-1}\alpha_i-A_{1\ldots i})=0.
$$
If $C_{i+1}\cap\cdots\cap C_k,\ldots, C_k, Y$ are all complex then the integrals over these cycles are all $0$, allowing us to replace $Y$ by $Y^*=C_{i+1}\cap\cdots\cap C_k$, and proving the last part of d), and thus (2.18).
\end{proof}
\section{Relation with Cheeger-Simons Differential Characters}

We recall that a differential character $\chi$ (of degree $r$) on a differentiable manifold $Y$ is a homomorphism from the group of differentiable singular $r$-cycles (with integer coefficients) to $\R/\Z$, such that there exists a differential $r+1$ form $\omega$ on $Y$ satisfying:  if $b$ is an $r$-boundary, say $b=\p c$, then
$$
\chi(b)=\int_{c}\omega\  (\textnormal{reduced mod }\Z).
$$
Further, there is a product operation on differential characters.

On a compact Riemannian manifold we can attach a differential character $\chi=\chi_A$ to an $n-r-1$ cycle $A$ in two steps:  first, we consider $r$-cycles $B$ whose support is disjoint from $A$ and define, using (1.2), 
$$
\chi(B)=(A,B)_Y\mod\Z.
$$
Second, following Cheeger's paper we look at any $B$, possibly intersecting $A$, and by transversality considerations find $B'$ homologous to $B$ and disjoint from $A$, and define
$$
\chi(B)=\chi(B')+\int_{D}\omega_A(\textnormal{mod }\Z)
$$
where $B'-B=\p D$, and $\omega_A$ is the harmonic form Poincare dual to the homology class of $A$.  It is easy to show that this definition does not depend on the choice of $B'$.

Now if we have several cycles $C_1,\ldots, C_k$ as in Section 1 and 2 we get corresponding characters $\chi_1,\ldots,\chi_k$.  The product of these characters is given by integrating (2.3) and taking limit as $t\to 0$ (except that the formulas in the papers on differential characters use different sign conventions).  Finally the product of the characters $\chi_1,\ldots,\chi_k$ is evaluated on $Y$, to obtain $(C_1,\ldots,C_k)\mod\Z$.

\end{document}